\documentclass[10pt]{amsart}
\usepackage{amssymb,amsmath,amsthm,graphicx,verbatim,amsbsy}

\newtheorem{theorem}{Theorem}[section]
\newtheorem{lemma}[theorem]{Lemma}

\theoremstyle{definition}

\theoremstyle{remark}
\newtheorem{remark}[theorem]{Remark}

\numberwithin{equation}{section}

\newcommand{\C}{\mathbb{C}}                                      %
\newcommand{\G}{{\Gamma}}                                        %
\newcommand{\Phit}{\tilde{\Phi}}                                   %
\newcommand{\cK}{\mathcal{K}}                                    %

\def\TT{{\mathbb T}}                                      %
\begin{document}


\title[Bivariate Bernstein-Szego measures on the
bi-circle]{Parameters associated with bivariate Bernstein-Szego measures on the
bi-circle}

\author{Jeffrey~S.~Geronimo}
\address{School of Mathematics, Georgia Institute of Technology, Atlanta, GA 30332-0160, USA}
\email{geronimo@math.gatech.edu}
\author{Philip Benge$^*$}\thanks{$^*$Author was an REU student during
  Summer 2010 suppported by grant DMS-0739343.}
\address{4214 Swire Ave., Apt 8, Baton Rouge, LA 70808, USA }
\email{pbenge2@tigers.lsu.edu}

\begin{abstract}
We consider measures supported on the bi-circle and review the recurrence relations satisfied by the orthogonal polynomials associated with these measures constructed using the lexicographical or reverse lexicographical ordering. New relations are derived among these recurrence coefficients. We extend the results of \cite{GW} on a parameterization for Bernstein-Szego measures supported on the bi-circle. 
\end{abstract}

\maketitle

\section{Introduction}

In this paper we continue the investigation begun in \cite{GW} on the orthogonal 
polynomials associated with measures supported on the bi-circle. In more than 
one variable an important consideration is which ordering to use. The usual 
ordering is the one suggested by Jackson \cite{J} which is the the total degree 
ordering. This is natural since the addition of any new polynomials does not 
alter the previous orthogonal polynomials already constructed. However in their 
solution of the two-variable Fejer-Reisz problem Geronimo and 
Woerdeman \cite{GW1}  were led to consider polynomials obtained using  the 
lexicographical or reverse lexicographical ordering (for an alternative viewpoint see  Knese \cite{K}). Orthogonal polynomials 
obtained using these orderings were first studied by Delsarte et al \cite{DGKa} who used them to solve the half-plane least squares problem \cite{DGKb}. 
Important in their work and later emphasized in \cite{GW} is the fact that in 
these orderings the moment matrices have a doubly Toeplitz structure. This 
allows a connection between the polynomials obtained using the above orderings 
and matrix orthogonal polynomials on the unit circle \cite{DGKa} 
(see also \cite{GW}). The various convergence properties of these polynomials in 
a strip as well as their connection to generalized Schur 
representations and Adamjan, Arov, and Krein theory were developed in \cite{DGKd}, and \cite{ DGKe}. 

Given a positive Borel probability measure $\sigma$ supported 
on the unit circle with an infinite number of points of increase let $\{\phi_n\}_{n\ge 0}$ be the sequence of polynomials of exact degree $n$ in $z=e^{i\theta}$ with positive leading coefficient having the property,
$$
\int_{\TT}\phi_j(e^{i\theta})\overline{\phi_k(e^{i\theta})}d\sigma(\theta)=\delta_{j,k}
$$
where $\delta_{j,k}$ is the Kronecker delta. These polynomials are known to
satisfy the recurrence formula \cite{G}, \cite{S}, \cite{Sz}, 
$$
\phi_{n}(z)=a_n(z\phi_{n-1}(z)-\alpha_n\overleftarrow\phi_{n-1}(z)),\ n\ge 1
$$
where $\overleftarrow\phi_n(z)=z^n\bar\phi_n(1/z)$ is called the reverse polynomial. The $\alpha_n$ are called the recurrence coefficients and $a_n=\frac{k_{n-1}}{k_n}$ where $k_n$ is the leading coefficient of $\phi_n$. From the orthogonality properties of $\phi_n$ and $\phi_{n-1}$ it is not difficult to obtain the relation
$$
1=a_n^2(1-|\alpha_n|^2).
$$
Thus $|\alpha_n|<1$ and given $\alpha_n$, $a_n$ can be computed. The recurrence 
coefficients play an 
important role in theory of orthogonal polynomials on the unit circle as can be 
seen by Verblunsky's Theorem \cite{S}.

\begin{theorem} Let $\sigma$ be a Borel probability supported on the unit circle with an infinite support, then associated with $\sigma$ is a unique sequence of recurrence coefficients $\{\alpha_n\}_{n=1}^{\infty}$ with $|\alpha_n|<1$. This correspondence is one-to-one.
\end{theorem}
A useful characterization theorem is the following,
\begin{theorem} Let $\sigma$ be a Borel  measure  supported on the unit circle. 
Then $\sigma$ is absolutely continuous with respect to Lebesgue measure with 
density $\frac{1}{|p_n(z)|^2}$ where $p_n(z)$ is a polynomial of exact degree $n$ in $z$ with $\overleftarrow p_n(z)$ nonzero for $|z|\le1$ if and only if $\alpha_i=0$, for $i>n$.
\end{theorem}

Measures of the form given in the above Theorem have come to be called
Bernstein-Szego measures \cite{S}.

In \cite{GW} a parameterization of the two variable trigonometric moment problem 
was introduced in an attempt to be able to find extensions of the above two 
Theorems to the two variable case. Here we continue the study of the algebraic 
properties of this problem. In section~2 we 
call together the results needed. In particular beginning with Borel measures 
supported on the bi-circle whose moment matrices are positive we construct 
orthogonal polynomials using the lexicographical or reverse lexicographical 
ordering. Then recurrence relations satisfied by these polynomials are 
displayed and some properties of the recurrence coefficients are noted. Here the 
parameterization discussed above is introduced and a two variable analog of 
Verblunsky's Theorem is presented. In section~3 we develop new equations 
between  the recurrence coefficients that shed light on how these coefficients 
are related to each other. Some of these recurrence relations are used in section~4 to develop an algorithm different from that given in \cite{GW} which allows us to make more precise the construction of the parameters left undetermined in Theorem 7.9 of \cite{GW}.  

\section{Preliminaries}
In this section we collect some results that will be used later. As noted above we will use the lexicographical ordering  which is defined by
$$
(k,\ell)<_{\rm lex} (k_1,\ell_1)\Leftrightarrow k<k_1\mbox{ or }
(k=k_1\mbox{ and } \ell<\ell_1),
$$
and the reverse lexicographical ordering, defined by
$$
(k,\ell)<_{\rm revlex} (k_1,\ell_1)\Leftrightarrow
(\ell,k)<_{\rm lex} (\ell_1,k_1).
$$
Both of these orderings are linear orders, and in addition they satisfy
$$
(k,\ell)<(m,n)\Rightarrow (k+p,\ell+q)<(m+p,n+q).
$$
In such a case, one may associate a half-space with the ordering
 which is defined by
 $\{ (k,l) \ : \  (0,0)  <  (k,l) \}$. In the case of the
 lexicographical ordering we shall denote the associated half-space by $H$
 and refer to it
 as {\it the standard half-space}. In the case of the reverse lexicographical
 ordering we shall denote the associated half-space by $\tilde H$. 
Let $\sigma$ be a positive Borel measure support on the 
bi-circle $z=e^{i\theta}$, $w=e^{i\phi}$ with Fourier coefficients,
$$
c_{k,j}=\int_{\TT} e^{-ik\theta} e^{-ij\phi}d\sigma(\theta,\phi).
$$
We now form the $(n+1)(m+1)\times(n+1)(m+1)$ moment matrix $C_{n,m}$ using the
lexicographical ordering. As noted in the introduction it has a special block Toeplitz form

\begin{equation}\label{toepone}
C_{n,m} = \left[
\begin{matrix}
C_0 & C_{-1} & \cdots & C_{-n}
\\
C_1 & C_0 & \cdots & C_{-n+1}
\\[-3pt]
\vdots &  & \ddots & \vdots \cr C_n & C_{n-1} & \cdots & C_{0}
\end{matrix}
\right],
\end{equation}
where each $C_i$ is an $(m+1)\times(m+1)$ Toeplitz matrix as
follows:
\begin{equation}\label{toeptwo}
C_i=\left[
\begin{matrix}
c_{i,0}& c_{i,-1} & \cdots & c_{i,-m}
\\
\vdots & &\ddots & \vdots \cr c_{i,m} &  & \cdots &
c_{i,0}
\end{matrix}
\right],\qquad i=-n,\dots, n.
\end{equation}
Thus $C_{n,m}$ has a doubly Toeplitz structure. If the reverse
lexicographical ordering is used in place of the lexicographical
ordering, we obtain another moment matrix $\tilde C_{n,m}$ where
the roles of $n$ and $m$ are interchanged. Throughout the rest of the
paper we will assume that $C_{n,m}$ is positive definite for all $0\le
n,\ m$.

We now compute orthogonal polynomials associated with $\sigma$. We begin by
ordering the monomials $z^i w^j,\ 0\le i\le n,\ 0\le j\le m$ lexicographically
then performing the Gram--Schmidt procedure using this ordering. Define the orthonormal polynomials
$\phi_{n,m}^l(z,w),\ 0\le n,\ 0\le m, \ 0\le l\le m,$
by the equations
\begin{equation}\label{sorthogonal}
\begin{split}
&\int_{\TT}\phi_{n,m}^l z^{-i}w^{-j}d\sigma=0, \quad  0\le i<n\ {\rm and }\
0\le j\le m\quad {\rm\ or}\ i=n\ {\rm and }\ 0\le j< l,\hspace*{-4pt}\\
&\int_{\TT}\phi_{n,m}^l\overline{\phi_{n,m}^l}=1,
\end{split}
\end{equation}
and
\begin{equation}\label{sorthogonaldeg}
\phi_{n,m}^{l}(z,w) = k^{n,l}_{n,m,l} z^n w^l + \sum_{(i,j)<_{\rm
lex}(n,l)} k^{i,j}_{n,m,l}z^iw^j.
\end{equation}
With the convention $k^{n,l}_{n,m,l}>0$, the above equations uniquely
specify $\phi^l_{n,m}$. Polynomials orthonormal with respect to
$\sigma$ but using the reverse lexicographical ordering will be
denoted by $\tilde \phi^l_{n,m}$. They are uniquely determined by the
above relations with the roles of $n$ and $m$ interchanged.

Set
\begin{equation}\label{vectpoly}
\Phi_{n,m}=\left[\begin{matrix} \phi_{n,m}^{m}\\ \phi_{n,m}^{m-1}\\[-2pt]
\vdots\\ \phi_{n,m}^{0} \end{matrix}\right] =
K_{n,m}\left[\begin{matrix} z^n w^m\\ z^n w^{m-1}\\[-2pt] \vdots\\ 1\end{matrix}\right],
\end{equation}
where the $(m+1)\times(n+1)(m+1)$ matrix $K_{n,m}$ is given by
\begin{equation}\label{knmmatrix}
K_{n,m}=\left[\begin{matrix} k_{n,m,m}^{n,m}&
k_{n,m,m}^{n,m-1}&\cdots & \cdots&\cdots& k_{n,m,m}^{0,0} \\ 0&
k_{n,m,m-1}^{n,m-1}&\cdots & \cdots&\cdots&
k_{n,m,m-1}^{0,0}\\\vdots &\ddots&\ddots&\ddots&\ddots&\ddots
\\0&\cdots& k_{n,m,0}^{n,0}& k_{n,m,0}^{n-1,m}&\cdots &
k_{n,m,0}^{0,0}\end{matrix}\right].
\end{equation}
As indicated above denote
\begin{equation}\label{vectpolytilde}
\tilde\Phi_{n,m} =\left[\begin{matrix} \tilde \phi_{n,m}^{n}\\ \tilde
\phi_{n,m}^{n-1}\\ \vdots\\ \tilde \phi_{n,m}^{0}\end{matrix}\right]
=\tilde K_{n,m}\left[\begin{matrix}w^m z^n \\  w^m z^{n-1}\\[-3pt] \vdots\\ 1
\end{matrix}\right],
\end{equation}
where the $(n+1)\times(n+1)(m+1)$ matrix $\tilde K_{n,m}$ is
given similarly to \eqref{knmmatrix} with the roles of $n$ and $m$
interchanged. For the bivariate polynomials $\phi^l_{n,m}(z,w)$ above
we define the reverse polynomials $\overleftarrow \phi^l_{n,m}(z,w)$ by the relation
\begin{equation}\label{reversepoly}
\overleftarrow \phi^l_{n,m}(z,w)=z^n w^m\bar{\phi}_{n,m}^l(1/z,1/w).
\end{equation}
With this definition $\overleftarrow \phi^l_{n,m}(z,w)$ is again a
polynomial in $z$ and $w$, and furthermore
\begin{equation}\label{reversevect}
\overleftarrow\Phi_{n,m}(z,w):= \left[\begin{matrix} \overleftarrow\phi_{n,m}^{m}\\ \overleftarrow\phi_{n,m}^{m-1}\\[-2pt]
\vdots\\ \overleftarrow\phi_{n,m}^{0} \end{matrix}\right]^T.
\end{equation}
An analogous procedure is used to define $\overleftarrow{\tilde\phi}^l_{n,m}$ and $\tilde\Phi_{n,m}$.

To find recurrence formulas for the vector
polynomials $\Phi_{n,m}$, we introduce the notation for every vector valued
polynomials $X$ and $Y$,
\begin{equation}\label{innerprod}
\langle X, Y\rangle=\int_{\TT} X(z,w) Y^{\dagger}(z,w)d\sigma, \quad |z|=1=|w|
\end{equation}
The following recurrence formulas which follow from the orthogonality
relations satisfied by $\Phi_{n,m}$ and $\Phit_{n,m}$ were proved in \cite{GW}.

\begin{theorem}\label{recurrencefor}
Given $\{\Phi_{n,m}\}$ and $\{\Phit_{n,m}\}$, $0<n$,
$0< m$, the following recurrence formulas hold:
\begin{align}
& A_{n,m}\Phi_{n,m} = z\Phi_{n-1,m} - \hat E_{n,m}\overleftarrow{\Phi}_{n-1,m}^T , \label{E1}
\\
&\Phi_{n,m}+ A^{\dagger}_{n,m}\hat E_{n,m}(A^T_{n,m})^{-1}\overleftarrow{\Phi}_{n,m}^T=
A^{\dagger}_{n,m}z\Phi_{n-1,m}, \label{E2}
\\
&\G_{n,m} \Phi_{n,m} = \Phi_{n,m-1} - {\cK}_{n,m} \Phit_{n-1,m},
\label{KK}
\\
&\G_{n,m}^1 \Phi_{n,m} = w \Phi_{n,m-1} - {\cK}^1_{n,m}\overleftarrow{\Phit}_{n-1,m}^T,
\label{K1}
\\
&\Phi_{n,m}=I_{n,m} \Phit_{n,m} + \G^{\dagger}_{n,m} \Phi_{n,m-1}
\label{II},
\\
&\overleftarrow{\Phi}_{n,m}^T=I_{n,m}^1 \Phit_{n,m} + (\G^1_{n,m})^T\overleftarrow{\Phi}_{n,m-1}^T
\label{II1},
\end{align}
where
\begin{align}
\hat E_{n,m} & = \langle z\Phi_{n-1,m}, \overleftarrow{\Phi}_{n-1,m}^T\rangle=\hat E_{n,m}^T \in
M^{m+1,m+1}, \label{enm}
\\
A_{n,m} & = \langle z\Phi_{n-1,m},\Phi_{n,m}\rangle \in
M^{m+1,m+1}, \label{anm}
\\
{\cK}_{n,m} & = \langle \Phi_{n,m-1}, \Phit_{n-1,m}\rangle \in
M^{m,n}, \label{knm}
\\
\G_{n,m} & = \langle \Phi_{n,m-1}, \Phi_{n,m} \rangle \in M^{m,m+1},
\label{gnm}
\\
{\cK}^1_{n,m} & = \langle w \Phi_{n,m-1}, \overleftarrow{\Phit}_{n-1,m}^T\rangle \in M^{m,n},
\label{k1nm}
\\
\G^1_{n,m} & = \langle w \Phi_{n,m-1}, \Phi_{n,m} \rangle \in
M^{m,m+1}, \label{g1nm}
\\
I_{n,m} & = \langle \Phi_{n,m}, \Phit_{n,m}\rangle \in
M^{m+1,n+1}, \label{inm}
\\
I^1_{n,m} & = \langle \overleftarrow{\Phi}_{n,m}^T, \Phit_{n,m}\rangle \in
M^{m+1,n+1}. \label{i1nm}
\end{align}
\end{theorem}\unskip

Here $M^{i,j}$ denotes the set of $i\times j$ matrices with complex entries. Equation~(\ref{E1}) was first found by Delsarte et. al. \cite{DGKa}.
\begin{remark}
{\rm Formulas similar to \eqref{E1}--\eqref{II1} hold for $\Phit_{n,m}$ 
and will be denoted by (\~{\ref{E1}})--(\~{\ref{II1}}). Throughout the rest of the 
paper we use the same notation to denote the tilde analogues of existing formulas stated for $\Phi_{n,m}$.}
\end{remark}

Examination of equation~(\ref{KK}) shows that the $(i,j)$ entries of $\Gamma_{n,m}$ are zero for $i\ge j$ with the $(i,i+1)$ entries positive. Likewise    equation~(\ref{K1}) implies  that entries $(i,j)$ of $\Gamma^1_{n,m}$ are zero for $i> j$ with the $(i,i)$ entries positive.

From the definitions of $\cK$, $\cK^1$, $I$, $I^1$ and their tilde analogs it is not difficult to see the following relations,
\begin{align}
 &{\tilde\cK}_{n,m}={\cK}_{n,m}^{\dagger},\ \tilde I_{n,m}=I_{n,m}^{\dagger},\label{3.24}
\\
& \tilde I^1_{n,m}=(I^1_{n,m})^T,\ \tilde{\cK}_{n,m}^1=({\cK}^1_{n,m})^T. \label{3.25}
\end{align}
Also the recurrence relations yield,
\begin{align}
&A_{n,m}A^{\dagger}_{n,m}=I_m-\hat E_{n,m}\hat E_{n,m}^{\dagger},\label{3.26}
\\
& \G_{n,m}\G_{n,m}^\dagger=I_m-{\cK}_{n,m} {\cK}_{n,m}^\dagger, \label{3.27}
\\
& \G^1_{n,m}(\G^1_{n,m})^\dagger=I_m- {\cK}^1_{n,m} ({\cK}^1_{n,m})^\dagger, \label{3.28}
\\
& I_{n,m}I^{\dagger}_{n,m}+\G^{\dagger}_{n,m}\G_{n,m}=I_{m+1}, \label{3.29}
\\
& I^1_{n,m}(I^1_{n,m})^\dagger+(\G^1_{n,m})^{\dagger}\G^1_{n,m}=I_{m+1}.
\label{3.30}
\end{align}

The definition of $\Phi_{n,m}$ implies that $A_{n,m}$ is an upper triangular 
matrix with positive diagonal entries. Thus it may be computed from 
$\hat E_{n,m}$ using a Cholesky decomposition of equation~(\ref{3.26}). A 
slightly more involved analysis \cite{GW} shows that $\Gamma_{n,m}$ may 
be computed using a Cholesky decomposition of  (\ref{3.27}).

The above recurrence formulas also give  pointwise formulas for
the recurrence coefficients. In order to obtain these formulas we
define the $m\times m+1$ matrices $U_m$ and $U^1_m$ as
\begin{equation}\label{um}
U_m=\left[ \begin{matrix}0,&I_m\end{matrix}\right],
\end{equation}
and
\begin{equation}\label{u1m}
U^1_m=\left[ \begin{matrix}I_m,&0 \end{matrix}\right],
\end{equation}
where $I_m$ is the $m\times m $ identity matrix. From equations~(\ref{matrixpolyz}) and (\ref{matrixpolyw}) we write
\begin{align}
\nonumber&\Phi_n^m(z) = \Phi^m_{n,n} z^n + \Phi^m_{n,n-1} z^{n-1} +\cdots,\\&\tilde\Phi_m^n(w) =
\tilde\Phi^n_{m,m} w^m + \tilde\Phi^n_{m,m-1} w^{m-1} +\cdots,\label{eqforonephi}
\end{align}
then the following relations hold:
\begin{align}
 &\G_{n,m}=\Phi^{m-1}_{n,n} U_m(\Phi^m_{n,n})^{-1},\label{3.32}
\\
& \G^1_{n,m}=\Phi^{m-1}_{n,n} U_m^1(\Phi_{n,n}^m)^{-1}, \label{3.33}
\\
& {\cK}_{n,m}=-\G_{n,m}I_{n,m}\tilde F_{n,m}, \label{3.34}
\\
& {\cK}^1_{n,m}=-\G^1_{n,m}\bar I^1_{n,m}\bar {\tilde F}^1_{n,m},
\label{3.35}
\end{align}
where $\tilde F_{n,m}=\tilde\Phi_{m,m}^n U_n^T (\tilde\Phi^{n-1}_{m,m})^{-1}$, and $\tilde
F^1_{n,m}=\tilde\Phi_{m,m}^n (U_n^1)^{T}(\tilde\Phi^{n-1}_{m, m})^{-1}$.

Equations~({\ref{E1}}) and ({\ref{E2}}) are a consequence of the relation between
$\Phi_{n,m}$ and the matrix orthogonal polynomials associated with the
$(m+1)\times(m+1)$ matrix measure $M_m$ given by
$$
dM_m(\theta)=\int_{\phi=-\pi}^{\pi}\left[\begin{matrix} w^m\\[-4pt] \vdots\\
1\end{matrix}\right]d\mu(\theta,\phi)\left[\begin{matrix} w^m\\[-4pt] \vdots\\
1\end{matrix}\right]^{\dagger},
$$
where $w=e^{i\phi}$. Given $E_{i,m}$ equation~({\ref{E1}}) allows the
computation of $\Phi_{i,m}$ along the strip $0\le i$, $0\le j\le m$. More precisely if we write
\begin{equation}\label{matrixpolyz}
\Phi_{n,m}(z,w)=\Phi_n^m(z)[w^m,\ldots,1]^T,
\end{equation}
then the $\Phi_i^m$ are a sequence of matrix polynomials of degree $i$ in $z$
satisfying
$$
\int_{-\pi}^{\pi}\Phi_i^m(z)dM(\theta)(\Phi_j^m(z))^{\dag}=I_{m+1}\delta_{i,j},
$$
where $I_{m+1}$ is the $(m+1)\times(m+1)$ identity matrix and $\delta_{i,j}$ is the Kronecker delta. Thus~({\ref{E1}}) and ({\ref{E2}}) follow the recurrence formulas satisfied by matrix polynomials orthogonal on the unit circle.  Equation~(\ref{vectpoly}) implies that the coefficient of
$z^i$ in $\Phi_i^m$, $\Phi^m_{i,i}$, is an $(m+1)\times(m+1)$ upper triangular matrix with 
positive diagonal entries. Similar statements hold for
\begin{equation}\label{matrixpolyw}
\Phit_{n,m}(z,w)=\Phit_n^m(w)[z^n,\ldots,1]^T.
\end{equation}
In contrast given $\Phi_{i,j}$,and $\Phit_{i,j}$ for $(i,j)=(n-1,m)$ or $(n,m-1)$
equations~({\ref{KK}}) and ({\ref{K1}}) allow the computation of
$\Phi_{n,m}$. 

As noted in \cite{GW} there is a lot of redundancy in the coefficients
of the above equations. If we have all the Fourier coefficients in the
notched rectangle $\{(i,j), 0\le i\le n,0\le j\le m\}\setminus(n,m)$ then
the polynomials $\Phi_{n,m-1},\  \tilde \Phi_{n-1 m}$ can be
computed. Notice that only two new Fourier
coefficients are required to compute all the polynomial $\Phi_{n,m}$  whereas $\cK_{n,m}$ and ${\cK}^1_{n,m}$ are both $m\times n $ matrices. This led in \cite{GW} to the introduction of parameters $u_{i,j}$
$i\ge0$, $u_{-i,-j}=\bar u_{i,j}$ such
that
\begin{equation}\label{konefirstent}
u_{-n,-m}=(\cK^1_{n,m})_{1,1}\ n>0,\ m>0 ,
\end{equation}
and
\begin{align}\label{defk}
u_{-n,m}&=(e^m_m)^T(\Phi^{m-1}_{n,n})^{-1}\cK_{n,m}((\tilde{\Phi}^{n-1}_{m,m})^{\dagger})^{-1}e^n_n\nonumber\\&=(e^m_m)^T\cK_{n,m}e^n_n/(k^{n,0}_{n,m-1,0}\tilde
k^{m,0}_{n-1,m,0}),\ n>0, \ m>0,
\end{align}
where $e^m_m$ is the $m$-dimensional vector with zeros in all its entries
except the last, which is one, and $k^{n,0}_{n,m-1,0}$ and $\tilde
 k^{m,0}_{n-1,m,0}$ are the leading coefficient of $\phi_{n,m-1}^0$ and
$\tilde\phi_{n-1,m}^0$ respectively. The last equality was
obtained using the upper triangularity of $\Phi^{m-1}_{n,n}$ and
$\tilde\Phi^{n-1}_{m,m}$, and equations~(\ref{vectpoly}) and
(\ref{vectpolytilde}). In terms of inner products the parameters can be
written as
\begin{equation}\label{innnerfircoeff}
u_{-n,-m}=\int_{\TT^2}\phi^{m-1}_{n,m-1}(z,w)\overline{\overleftarrow{\tilde\phi^{n-1}_{n-1,m}}(z,w)}d\sigma(\theta,\phi),\
z=e^{i\theta},\ w=e^{i\phi},
\end{equation}
and
\begin{equation}\label{innerseccoeff}
u_{-n,m}=\int_{\TT^2}\hat\phi^{0}_{n,m-1}(z,w)\overline{\hat{\tilde\phi}^{0}_{n-1,m}(z,w)}d\sigma(\theta,\phi),\
z=e^{i\theta},\ w=e^{i\phi}.
\end{equation}
Here $\hat\phi^{0}_{n,m-1}$ and $\hat{\tilde\phi}^{0}_{n,m-1}$ have leading
coefficient one. Since $\cK_{n,m}$ and $\cK^1_{n,m}$ are contractions the
parameters must satisfy the constraints 
$$|u_{n,m}|<1$$ 
and 
$$k^{n,0}_{n,m-1,0}\tilde
k^{m,0}_{n-1,m,0}|u_{n,-m}|<1.$$

With this the following Theorem was proved
in \cite{GW},

\begin{theorem}\label{th6.1} Given parameters
$u_{i,j}\in\C$, $0\le i$, $u_{-i,j}=\bar u_{i,-j}$ construct
\begin{itemize}
\item  scalars $\hat E_{i,0}$ and
$\tilde{\hat E}_{0,j}$;
\item matrices $\cK_{i,j}$, $i>0$, $j>0$; and 
\item  numbers $(e^j_1)^T H^3_{i,j}e^j_1$, $i>0$, $j>0$. 
\end{itemize}
If
\begin{equation}\label{6.1}
u_{0,0}>0,\ |\hat E_{i,0}|<1, |\tilde{\hat E}_{0,j}|<1,\ ||\cK_{i,j}||<1, \text{ and } (H^3_{i,j})_{1,1}<1,
\end{equation}
then there exists a unique positive measure $\sigma$ supported on the bi-circle such that
\begin{equation}\label{6.2}
\int_{\TT}\Phi_{i,m}\Phi_{j,m}^{\dagger}d\sigma=\delta_{i,j}I_{m+1}
\text{ and }
\int_{\TT}\Phit_{n,i}\Phit_{n,j}^{\dagger}d\sigma=\delta_{i,j}I_{n+1}.
\end{equation}
The conditions \eqref{6.1} are also necessary.
\end{theorem}
The numbers $(H^3_{i,j})_{1,1}$ are given by equation~(5.14) in \cite{GW} .

This Theorem is the two dimensional analog of Verblunsky's Theorem
discussed in the introduction. A polynomial $p$ is of degree $(n,m)$
if
$$
p(z,w)=\sum_{i=0}^n\sum_{j=0}^m p_{i,j} z^i w^j,
$$
with $p_{n,m}\ne0$.

When $d\sigma=\frac{1}{|p_{n,m}(z,w)|^2}$ where
$p_{n,m}$ is a polynomial of degree $(n,m)$
with $\overleftarrow p_{n,m}$ stable (i.e.\ $\overleftarrow p_{n,m}\ne0$, $|z|,|w|\le 1$) then more can be
said.
\begin{theorem}\label{verbtwod} Let $\mu$ be a positive measure on the 
bicircle. Then $\mu$ is purely absolutely continuous with respect to the  
Lebesgue measure and $d\mu=\frac{d\theta d\phi}{4\pi^2 |p_{n,m}|^2}$, where
$p_{n,m}$ is a polynomial of degree $(n,m)$ with $\overleftarrow p_{n,m}$ 
stable if and only if

{\rm (a)} $\cK_{n,j}=0$, $\tilde{\hat E}_{n-1,j+1}=0$, and $u_{n,j+1}=0$, $j\ge m$; 

{\rm (b)} $\cK_{i,m}=0$, $\hat E_{i,m-1}=0$, and $u_{i,m}=0$, $i>n$; 

{\rm (c)} $u_{i,j}=0,|i|>n,\ j>m$.
\end{theorem}
It was actually shown that in this case $u_{-n,j}$, $u_{-i,m}$,
$u_{n-1,j+1}$ and $u_{i,m-1}$ are equal to zero for $j\ge m,\ i>n$. How to
compute the remaining parameters and how they were related to $u_{i,j}$, where
$|i|\le n$ and $j\le m$, was not indicated and it the subject of the
remaining sections.

\section{Relations for $\hat E_{n,m}$, $\cK_{n,m}$ and ${\cK}^1_{n,m}$}
In order to prove Theorem~\ref{th6.1} it was necessary to show that most of
the entries in $\cK_{n,m}$ and ${\cK}^1_{n,m}$ could be computed knowing
the recurrence coefficients on the $(n-1,m)$ and $(n, m-1)$ levels. These relations
will be augmented by the following new relations which will be used later to compute coefficients on lower levels from those on higher levels.

\begin{lemma}\label{enma} For $n>0$ and $m\ge0$,
\begin{equation}\label{eenm} 
\hat E_{n,m}=\Gamma_{n-1,m+1}\hat E_{n,m+1}(\Gamma^1_{n-1,m+1})^T+{\cK}_{n-1,m+1}({\cK}^1_{n-1,m+1})^T.
\end{equation}
Also for $n\ge0$ and $m>0$,
\begin{equation}\label{teenma}
\hat{\tilde E}_{n,m}=\tilde\Gamma_{n+1,m-1} \hat{\tilde E}_{n+1,m}(\tilde\Gamma^1_{n+1,m-1})^T+\tilde{\cK}_{n+1,m-1}(\tilde{\cK}^1_{n+1,m-1})^T.
\end{equation}
\end{lemma}
\begin{proof}
 From the definition of $\hat E_{n,m}$,
\begin{equation}
\label{eequation}
\hat E_{n,m}=\langle z\Phi_{n-1,m},\overleftarrow\Phi_{n-1,m}^T\rangle
\end{equation}
eliminate $\Phi_{n-1,m}$ using equation~(\ref{KK}) to obtain
\begin{align*}
\hat E_{n,m}&=\Gamma_{n-1,m+1}\langle z\Phi_{n-1,m+1},\overleftarrow\Phi_{n-1,m}^T\rangle\\
&\quad +K_{n-1,m+1}\langle z\tilde\Phi_{n-2,m+1},\overleftarrow\Phi_{n-1,m}^T\rangle.
\end{align*}
With the use of the reverse of (\ref{K1}) the first integral on the right hand side of the above equation can be rewritten as
\begin{align*}
\langle z\Phi_{n-1,m+1},\overleftarrow\Phi_{n-1,m}^T\rangle
&=\langle z\Phi_{n-1,m+1},\overleftarrow\Phi_{n-1,m+1}^T\rangle
(\Gamma^1_{n-1,m})^T\\
&\quad +\langle\Phi_{n-1,m+1}, \Phi_{n-2,m+1}\rangle(\cK^1_{n-1,m+1})^T\\
&=\hat E_{n,m+1}(\Gamma^1_{n-1.m+1})^T.
\end{align*}
Equation~(\ref{enm})  and the orthogonality of $\Phi_{n-1,m+1}$ to $\tilde\Phi_{n-n,m+1}$ has been used to obtain the last equality. The result now follows by taking the transpose of equation~(\ref{k1nm}). Equation~(\ref{teenma}) follows in a similar manner using the  tilde analog of the above equations.
\end{proof}

\begin{lemma}\label{knma}
\begin{equation}\label{Kequation}
\Gamma_{n-1,m}\hat E_{n,m}I_{n-1,m}^{1}=A_{n,m-1}\cK_{n,m}-\cK_{n-1,m}\tilde{\Gamma}_{n-1,m}^{1}
\end{equation}
and
\begin{equation}\label{Ktequation}
\tilde\Gamma_{n,m-1}\hat{\tilde E}_{n,m}\tilde I_{n,m-1}^{1}=\tilde A_{n-1,m}\tilde \cK_{n,m}-\tilde \cK_{n,m-1}\Gamma_{n,m-1}^{1}
\end{equation}
\end{lemma}
\begin{proof}
To obtain (\ref{Kequation}) note that  equations~(\ref{enm}) and
the reverse transpose of (\~{\ref{II1}}) give,
$$
\hat E_{n,m}I_{n-1,m}^{1}=\langle z\Phi_{n-1,m},\tilde\Phi_{n-1,m}\rangle.
$$
Multiplying the above equation on the left by $\Gamma_{n-1,m}$ then using
the reverse transpose of equation~(\ref{KK}) yields,
$$
\Gamma_{n-1,m}\hat E_{n,m}I_{n-1,m}^{1}=
\langle z\Phi_{n-1,m-1},\tilde\Phi_{n-1,m}
\rangle -\cK_{n-1,m}\langle z\tilde\Phi_{n-2,m},\tilde\Phi_{n-1,m}\rangle.
$$  
The second integral in the above equation evaluates to $\tilde\Gamma_{n-1,m}^1$. Substitution of equation~(\ref{E1}) in the first integral to eliminate $z\Phi_{n-1,m-1}$ then using equation~(\ref{knm}) yields equation~(\ref{Kequation}). The argument for equation~(\ref{Ktequation}) follows in an analogous manner using the tilde analog of the above equations.
\end{proof}

Finally
\begin{lemma}\label{k1nma}
\begin{equation}\label{K1equation}
I_{n-1,m}^{\dag}\hat E_{n,m}(\Gamma_{n-1,m}^{1})^{T}=(\cK_{n,m}^{1})^{T}A_{n,m-1}^{T}-\tilde{\Gamma}_{n-1,m}^{\dagger}(\cK_{n-1,m}^{1})^{T}
\end{equation}
and
\begin{equation}\label{K1tequation}
\tilde I_{n,m-1}^{\dag}\hat{\tilde E}_{n,m}(\tilde\Gamma_{n,m-1}^{1})^{T}=(\tilde{\cK}_{n,m}^{1})^{T}\tilde A_{n-1,m}^{T}-\Gamma_{n,m-1}^{\dagger}(\tilde{\cK}_{n,m-1}^{1})^{T}
\end{equation}
\end{lemma}
\begin{proof}
To obtain (\ref{Kequation}) use   equations~(\ref{enm}) and (\~{\ref{II}}) to find
$$
 I^{\dag}_{n-1,m}\hat E_{n,m}=\langle z\tilde\Phi_{n-1,m},
\overleftarrow\Phi_{n-1,m}^T\rangle.
$$
Multiplying the above equation on the left by the transpose of $\Gamma^1_{n-1,m}$ then using
the reverse transpose of equation~(\ref{K1}) yields
\begin{align*}
I^{\dag}_{n-1,m}\hat E_{n,m}(\Gamma^1_{n-1,m})^{T}
&=\langle z\tilde\Phi_{n-1,m},\overleftarrow{\Phi}_{n-1,m-1}\rangle
\nonumber\\
&\quad -\langle \tilde\Phi_{n-1,m},\tilde\Phi_{n-2,m}\rangle
(\cK^1_{n-1,m})^T.
\end{align*}
The second integral in the above equation evaluates to $\tilde\Gamma^{\dag}_{n-1,m}$. Substitution of the reverse transpose of equation~(\ref{E1})  in the first integral then using equation~(\~{\ref{k1nm}}) yields equation~(\ref{K1equation}).
As above equation~(\ref{K1tequation}) follows a similar argument using the tilde analogs of the above equations.
\end{proof}

With these recurrences we can prove a strengthening of Lemma~7.5 in \cite{GW}

\begin{lemma}\label{eij} If $\hat E_{i,j}=0$, then the first column of $\cK^1_{i,j}$ is equal to zero, in particular $u_{i,j}=0$. If $\hat E_{i,j}$ and  $\cK_{i-1,j}(\cK^1_{i-1,j})^T$ are zero, then so is $\hat E_{i,j-1}$. If $\cK_{i,j},\ \hat E_{i,j-1}$, and $u_{i,j}$ are zero, then $\hat E_{i,j}=0$. 
Likewise if  $\hat {\tilde E}_{i,j}=0$, then the first row of $K^1_{i,j}$ is equal to zero. If $\hat{\tilde E}_{i,j}$, and $\cK_{i,j-1}^{\dag}\cK^1_{i,j-1}$ are zero, then so is $\hat{\tilde E}_{i,j-1}$. If $\cK_{i,j},\ \hat{\tilde E}_{i-1,j}$, and $u_{i,j}$ are zero, then $\hat{\tilde E}_{i,j}=0$.
\end{lemma}\unskip

\begin{proof} If $\hat E_{i,j}=0$, then (\ref{K1equation}) and the triangular 
structure of $\tilde \G_{i-1,j}$ show that the first column of $\cK^1_{i,j}$ is 
zero. If $\hat E_{i,j}$ and $\cK_{i-1,j}(\cK^1_{i-1,j})^T$ are equal to zero, then (\ref{eenm}) shows that $\hat E_{i, j-1}=0$. Equations~(3.50) and (3.51) in \cite{GW} are
\begin{align}
&\G_{n-1,m}\hat E_{n,m}=A_{n,m-1}\cK_{n,m}(I^1_{n-1,m})^{\dagger}+\hat E_{n,m-1}\bar\G^1_{n-1,m},
\label{EE1}
\\
&\hat E_{n,m}(\G^1_{n-1,m})^T=I_{n-1,m}(\cK^1_{n,m})^TA^T_{n,m-1}+\G^{\dagger}_{n-1,m}\hat
E_{n,m-1}.\label{EE2}
\end{align}
Thus if $\cK_{i,j}$ and $\hat{\tilde E}_{i-1,j}$ are equal to zero (\ref{EE1}) and the fact that
$\hat E_{i,j}$ is symmetric shows that all its entries are equal to zero except for $[\hat E_{i,j}]_{(1,1)}$. Equation~(\ref{EE2}) and the assumption that $u_{i,j}=0$ give that this entry is equal to zero also. The remaining statements 
follow in an analogous fashion using the tilde analogs of the above equations.

\end{proof}
For the next lemma we recast equations~(\ref{vectpoly}) and (\ref{vectpolytilde}) as,
\begin{equation}\label{lnmm}
\Phi_{n,m}(z,w)=\sum_{i=0}^m L_{n,m}^i w^i[z^n,\ldots,1]^T
\end{equation}
and
\begin{equation}\label{tlnmm}
\Phit_{n,m}(z,w)=\sum_{i=0}^n\tilde L_{n,m}^i z^i[w^m,\ldots,1]^T
\end{equation}
From this we have,
\begin{lemma}\label{pointk} For $n,m\ge0$
\begin{equation}\label{cominm}
I_{n,m}=L_{n,m}^m(\Phit_{m,m}^n)^{-1}
\end{equation}
and
\begin{equation}\label{comi1nm}
I^1_{n,m}=\bar L_{n,m}^0 J_n (\Phit_{m,m}^n)^{-1},
\end{equation}
where $J_n$ is the $(n+1)\times(n+1)$ matrix with ones on the reverse-diagonal
and zeros everywhere else.
\end{lemma}
\begin{proof}
The first formula follows from equation~(\ref{lnmm}). The second equation can be seen from the computation,
$$
\overleftarrow{\Phi}_{n,m}(z,w)^T=\sum_{i=0}^m \bar L_{n,m}^i w^{m-i}J_n[z^n,\ldots,1]^T.
$$

\end{proof}

The above results give formulas for the parameters. From equations~(\ref{3.32}), (\ref{3.34}), and (\ref{cominm}) we find
\begin{equation}\label{firpara}
(e^j_j)^T(\Phi^{j-1}_{i,i})^{-1}\cK_{i,j}((\Phit^{i-1}_{j,j})^{\dag})^{-1}e^i_i=-(e^{j+1}_{j+1})^T(\Phi^j_{i,i})^{-1}L^j_{i,j}U_i^T((\Phit^{i-1}_{j,j})^{\dag}\Phit^{i-1}_{j,j})^{-1}e^i_i,
\end{equation}
where we have used the fact the $U_j^T e^j_j=e^{j+1}_{j+1}$ in the last 
equation. Likewise
equations~(\ref{3.33}), (\ref{3.35}), and  (\ref{comi1nm})
show,
\begin{equation}\label{secpara}
(\cK^1_{i,j})_{1,1}=-(\Phi^{j-1}_{i,i}U_i(\Phi^j_{i,i})^{-1}\bar L_{i,j}^0 J_{i} (U_i^1)^T((\Phit^{i-1}_{j,j})^{\dag}\Phit^{i-1}_{j,j})^{-1})_{1,1}.
\end{equation}

\section{Construction of the Parameters}

Using the results above we are now able to compute the remaining parameters
from those given in the rectangle $0\le i\le n,\ 0\le j\le m$.
We first show that all the polynomials $\Phi_{i,j}$ and $\Phit_{i,j}$ can
be computed for $i>n$ and $0\le j\le m$.
To see this suppose that we are given $\Phi_{n,j}$ and $\Phit_{n,j}$ 
for  $0\le j\le m$ and also 
conditions (a), (b), and (c) of Theorem~\ref{verbtwod} are satisfied. From 
their defining properties we see that $\phi^m_{n,m}=\tilde\phi^n_{n,m}$. 
Lemma~\ref{eij} shows that
(a) and (b) imply that $\hat E_{i,m}=0$ for $i>n$ so that from 
equation~(\ref{E1}) $\Phi_{i,m}=z^{i-n}\Phi_{n,m}$ and because
$\hat E_{i,m-1}=0,\ i> n$, $\Phi_{i,m-1}=z^{i-n}\Phi_{n,m-1}$. If $\cK_{i,m}=0$ it 
follows from the triangularity of $\Gamma_{i,m}$ and equation~(\ref{3.27}) 
that $\Gamma_{i,m}=U_m$. Likewise $\tilde\Gamma_{i,m}=U_i$. This implies 
through equations~(\~{\ref{KK}}) and (\~{\ref{inm}}) that 
$$
\tilde \Phi_{i,m}=\begin{bmatrix} z^{i-n}\phi^m_{n,m}\\ \tilde\Phi_{i-1,m}
\end{bmatrix}
$$
for $i\ge n$ which gives all of $\tilde\Phi_{i,m}$ for $i> n$. 
Set $w=0$ in (\~{\ref{E2}}) then utilize equation~(\ref{matrixpolyw}) and the fact  that $\overleftarrow\Phit_i^n(0)$ is invertible to obtain
\begin{equation}\label{eae}
\tilde A^{\dag}_{i,m}\tilde{ \hat E}_{i,m}(\tilde A^T_{i,m})^{-1} = -\Phit^i_m(0)
 J_i(\overleftarrow{\Phit}^i_m(0)^{-1})^T\equiv B^i_m.
\end{equation}
To find $\tilde A_{i,m}$ use the orthogonality properties of $\Phit_{i,m}$
and $\overleftarrow\Phit_{i,m}$ in equation (\~{\ref{E2}}) to find,
\begin{equation*}
I-B^i_m(B^i_m)^{\dag}=\tilde A_{i,m}^{\dag}\tilde A_{i,m}.
\end{equation*}
Since $\tilde A_{i,m}$ is upper triangular with positive diagonal entries
it may be computed using the lower Cholesky factorization of the left hand
side of the above equation.  Using this in (\~{\ref{E2}}) allows us to
compute $\Phit_{i,m-1}$ for $i>n$. In an analogous fashion $\Phit_{i,j}$
may be computed for $i>n$ and $0\le j\le m-1$. Now $\tilde \Gamma_{i,j}$
and $\tilde\Gamma^1_{i,j}$, $i>n$, $0<j\le m-1$ may be computed from from
equations~(\~{\ref{3.32}}) and (\~{\ref{3.33}}) respectively.  With $i=m-1$
we find from (\ref{eenm}) since $\hat E_{n+1,m-1}=0$ that,
$$
\hat E_{n+1,m-2}=\cK_{n,m-1}(\cK^1_{n,m-1})^T
$$
which gives $\hat E_{n+1,m-2}$ because by assumption $\cK_{n+1,m-1}$ and $\cK^1_{n+1,m-1}$ are known. Since $A_{n+1,m-2}$ may be computed from $\hat E_{n+1,m-2}$ using the upper Cholesky factorization of (\ref{3.26}) we obtain $\Phi_{n+1,m-1}$ from (\ref{E1}). By induction we see that the above argument gives $\hat E_{n+1,i}$, $i=0,\ldots, m-3$ from which $A_{n+1,i}$ may be computed and then $\Phi_{n+1,i}$. Using equations~(\ref{cominm}), (\ref{comi1nm}), 
(\ref{3.32})--(\ref{3.35}) allows
us to compute $\Gamma_{n+1,i},\Gamma^1_{n+1,i}$, $I_{n+1,i}$,
$I^1_{n+1,i}$, $\cK_{n+1,i}$, and $\cK^1_{n+1,i}$ for $0<i\le m-1$. With
equation~(\ref{eenm}) and the coefficients just computed we repeat the above argument for level $(n+2,i)$ and by induction $(i,j)$, $i>n$, $0\le j\le m-1$. 

We summarize this with
\begin{lemma}\label{comcoeff} Given (a), (b), and (c) of 
Theorem~{\rm\ref{verbtwod}} as well as $\Phi_{n,j}$ and $\Phit_{n,j}$ for $0\le j\le m$ then $\Phi_{i,j}$ and $\Phit_{i,j}$ for $i>n, \ 0\le j\le m$ can be computed recursively. If $\Phi_{i,m}$ and $\Phit_{i,m}$ are given then $\Phi_{i,j}$ and $\Phit_{i,j}$ for $0\le i\le n, \ j> m$ can be computed recursively.
\end{lemma}

We now use the formulas (\ref{firpara}) and (\ref{secpara}) which give $u_{-i,j}$ for $i>n$ and $1\le j<m$ and  $u_{i,j}$ for $i>n$ and
$1\le j<m-1$. For $j=0$ we have from \cite{GW} that 
\begin{equation}\label{zpara}
u_{i,0}=-\frac{\Phi_{i,0}^0}{\Phi_{i,i}^0},
\end{equation} 
which gives the parameters in the strip $i>n,\ 0\le j\le m-1$.
To compute the parameters in the strip $j>m,\ 0\le i\le n-1$
equations~(\ref{3.24}) and (\ref{3.25}) show that we need only interchange
$i$ with $j$ and the matrices associated with the lexicographical ordering with
those associated with the reverse   lexicographical ordering in
equations~(\ref{firpara}), (\ref{secpara}), and (\ref{zpara}). 
This leads to,

\begin{theorem} Suppose $\sigma$ is a positive Borel measure supported on
  the bi-circle with parameters $u_{i,j}$. If

{\rm (a)} $u_{i,-j}, 1\le i, j\le n$ and  $u_{i,j}, 0\le i, j\le n$ give $\cK_{n,m}=0$;

{\rm (b)} $u_{i,j}=0$ for $i=n-1, j>m$, $i> n, j=m-1$, $|i|> n, j=m$, $|i|\ge n,j>m$;

{\rm (c)}  for $i>n$, $0\le j\le m-2$, $u_{i,j}$ are equal to the left hand sides of (\ref{firpara}), (\ref{secpara}), or (\ref{zpara})computed  using the above algorithm;

{\rm (d)} for $j>m$, $0\le i\le n-2$,  $u_{i,j}$ are equal to the left hand sides of (\~{\ref{firpara}}), (\~{\ref{secpara}}), of (\~{\ref{zpara}}) computed using the above algorithm,

\medskip

\noindent then $\sigma$ is absolutely continuous with respect to Lebesgue measure
with density $\frac{1}{|p_{n,m}|^2}$ where $p_{n,m}$ is of degree $(n,m)$
with $\overleftarrow p_{n,m}(z,w)$ stable. $p_{n,m}$ is unique up to
multiplication by a complex number of modulus one.
\end{theorem}  
\begin{proof} If $\cK_{n,m}=0$ then from Theorem~(7.3) of \cite{GW} or Theorem~(10.1) of \cite{K}
there exits a polynomial $p_{n,m}$ of degree $(n,m)$ with $\overleftarrow p_{n,m}$ stable such that the measure $d\rho=\frac{d\theta d\phi}{4\pi^2 |p_{n,m}|^2}$ has parameters given in (a). By Theorem~(\ref{verbtwod}) and the algorithm above we that this measure has the same paramters as in (b) and (c). Thus by Theorem~(\ref{th6.1}) $d\sigma= \frac{d\theta d\phi}{4\pi^2|p_{n,m}|^2}$.

\end{proof}

\end{document}